\documentclass[a4paper,11pt]{article}
\usepackage[english]{babel}
\usepackage[latin1]{inputenc}
\usepackage{amssymb}
\usepackage{amsmath}
\usepackage{amscd}
\usepackage{latexsym}
\usepackage{amsthm}
\usepackage{eucal}
\usepackage{eufrak}
\usepackage{psfig}
\theoremstyle{plain}
\newtheorem{thm}{Theorem}
\newtheorem{cor}[thm]{Corollary}
\newtheorem{lem}[thm]{Lemma}
\newtheorem{prop}[thm]{Proposition}

\theoremstyle{definition}
\newtheorem{defn}{Definition}
\theoremstyle{remark}
\newtheorem{oss}{Remark}
\newtheorem{esem}{Example}

\begin{document}

\title{\textbf{A Gabriel
Theorem for Coherent Twisted Sheaves}}

\author{Arvid Perego}

\maketitle

\begin{abstract}
The aim of this work is to give a generalization of Gabriel's
theorem for twisted sheaves over smooth varieties. We start by
showing that we can reconstruct a variety $X$ from the category
$Coh(X,\alpha)$ of coherent $\alpha-$twisted sheaves over $X$.
This follows from the bijective correspondence between closed
subsets of $X$ and Serre subcategories of finite type of
$Coh(X,\alpha)$. Then we show that any equivalence between
$Coh(X,\alpha)$ and $Coh(Y,\beta)$, where $X$ and $Y$ are smooth
varieties, induces an isomorphism between $X$ and $Y$. Here, the
problem is to show that we can extend any coherent twisted sheaf
on an open subscheme of $X$ to a coherent twisted sheaf on $X$. In
order to do this, we study perfect and compact objects in
$D(QCoh(X,\alpha))$. As a complement, we study the problem of
saturatedness of $D^{b}(X,\alpha)$, which will be proved at least
for smooth and proper varieties.
\end{abstract}

\section[Introduction]{Introduction}

Gabriel's theorem is one of the main results of the use of
category theory in algebraic geometry. It says that for every
noetherian scheme $X$, we can construct a scheme $E_{X}$ from
$Coh(X)$, and an isomorphism between $E_{X}$ and $X$, so that we
can say that $Coh(X)$ carries informations about the scheme
structure of $X$. Moreover, we have that two noetherian schemes
have equivalent categories of coherent sheaves if and only if they
are isomorphic. \indent What we want to do is to show that we can
extend this theorem to the case of $Coh(X,\alpha)$, the category
of coherent $\alpha-$twisted sheaves over $X$, where such $X$ is a
(smooth) $k-$variety, that is, a separeted scheme of finite type
over a field $k$. More precisely, we want to show the following:

\begin{thm}
\label{thm:gab}Let $X$ be a variety over a field $k$, and
$\alpha\in BrX$. Then the abelian category $Coh(X,\alpha)$
determines $X$. Moreover, if $X$ and $Y$ are two smooth varieties,
$\alpha\in BrX$, $\beta\in BrY$, any equivalence between
$Coh(X,\alpha)$ and $Coh(Y,\beta)$ induces an isomorphism between
$X$ and $Y$.
\end{thm}

Note that this theorem does not tell anything about how $f^{*}$
acts on $BrY$, namely if $f^{*}\beta=\alpha$. This would be an
easy consequence of the twisted version of Orlov's theorem which
is shown in \cite{CS} (see \cite{CS}, Remark 5.4 for the proof of
this fact).

In this introduction, we would like to recall the notions of
twisted sheaf and of category of (quasi) coherent twisted sheaves.
The main reference for definitions and proofs will be \cite{Cal}.
In the following, let $X$ be a variety over a field $k$. We will
denote by $Br'X:=H^{2}(X,\mathcal{O}^{*}_{X})_{tors}$ the
cohomological Brauer group of $X$, and by $BrX$ the Brauer group
of $X$, that is, the group of equivalence classes of Azumaya
algebras over $X$. From Theorem 1.1.8 in \cite{Cal}, we know that
$BrX$ is a subgroup of $Br'X$, so that, in particular, we will
think an element $\alpha\in BrX$ as the cohomology class of a
2-cocycle $\{\alpha_{ijk}\}\in
\check{C}^{2}(X,\mathcal{U},\mathcal{O}^{*}_{X})$, where
$\mathcal{U}=\{U_{i}\}_{i\in I}$ is an open covering of $X$. We
will denote $U_{ij}=U_{i}\cap U_{j}$ and $U_{ijk}=U_{i}\cap
U_{j}\cap U_{k}$, so that $\alpha_{ijk}\in
\Gamma(U_{ijk},\mathcal{O}^{*}_{X})$ and the 2-cocycle condition
is satisfied. In this way we have the following definition.

\begin{defn}
\label{defn:twisted}We call \textit{sheaf twisted by}
$\alpha\in\check{C}^{2}(X,\mathcal{U},\mathcal{O}^{*}_{X})$, or
simply $\alpha-$\textit{sheaf}, a family
$\mathcal{F}=(\mathcal{F}_{i},\varphi_{ij})_{i,j\in I}$ where
$\mathcal{F}_{i}$ is an $\mathcal{O}_{U_{i}}-$module and
$$\varphi_{ij}:\mathcal{F}_{j|U_{ij}}\stackrel{\sim}\longrightarrow\mathcal{F}_{i|U_{ij}}$$is
an isomorphism of $\mathcal{O}_{U_{ij}}-$modules such that:
\begin{enumerate}
\item $\varphi_{ii}=id_{\mathcal{F}_{i}}$ for every $i\in I$;
\item $\varphi_{ij}=\varphi_{ji}^{-1}$ for every $i,j\in I$; \item
$\varphi_{ij}\circ\varphi_{jk}\circ\varphi_{ki}=\alpha_{ijk}\cdot
id_{\mathcal{F}_{i|U_{ijk}}}$.
\end{enumerate}
If the sheaves $\mathcal{F}_{i}$ are quasi-coherent (coherent,
locally free) $\mathcal{O}_{U_{i}}-$modules for every $i\in I$, we
say that $\mathcal{F}$ is a \textit{quasi-coherent}
(\textit{coherent}, \textit{locally free})
$\alpha-$\textit{sheaf}.
\end{defn}

\begin{defn}
A \textit{morphism} between two $\alpha-$sheaves
$\mathcal{F}=(\mathcal{F}_{i},\varphi_{ij})$ and
$\mathcal{G}=(\mathcal{G}_{i},\psi_{ij})$ is a family
$f=(f_{i})_{i\in I}$ where
$$f_{i}:\mathcal{F}_{i}\longrightarrow\mathcal{G}_{i}$$is a
morphism of $\mathcal{O}_{U_{i}}-$modules such that
$\psi_{ij}\circ f_{j}=f_{i}\circ\varphi_{ij}$ for every $i,j\in
I$.
\end{defn}

We will write $Mod(X,\mathcal{U},\alpha)$ the category whose
objects are $\alpha-$sheaves over $X$ and morphisms are morphisms
of $\alpha-$sheaves. We define also its full subcategories
$QCoh(X,\mathcal{U},\alpha)$ and $Coh(X,\mathcal{U},\alpha)$ in
the obvious way.

Actually, we can show that if we change the open covering
$\mathcal{U}$ into $\mathcal{U}'$, we have that
$Mod(X,\mathcal{U},\alpha)$ is canonically equivalent to
$Mod(X,\mathcal{U}',\alpha)$, so that we are allowed to write
$Mod(X,\alpha)$ instead of $Mod(X,\mathcal{U},\alpha)$ (see
\cite{Cal}, Lemma 1.2.3 and Corollary 1.2.6). Moreover, we can
show that if we change the 2-cocycle $\alpha$ to an equivalent one
$\beta$, we can find a (non-canonical) equivalence between
$Mod(X,\alpha)$ and $Mod(X,\beta)$, so that $Mod(X,\alpha)$
depends only on $\alpha\in Br'X$ (see \cite{Cal}, Lemma 1.2.8).

It is easy to show that $Mod(X,\alpha)$, $QCoh(X,\alpha)$ and
$Coh(X,\alpha)$ are abelian categories. Moreover, $Mod(X,\alpha)$
and $QCoh(X,\alpha)$ have enough injective objects (see
\cite{Cal}, Lemma 2.1.1).

\begin{oss} If $\mathcal{F}=(\mathcal{F}_{i},\varphi_{ij})$ is an
$\alpha-$sheaf and $x\in U_{i}$, we write
$\mathcal{F}_{x}:=\mathcal{F}_{i,x}.$ Actually, this does not
define a true stalk for the twisted sheaf, since it is determined
only up to isomorphism. Anyway, the following definition makes
sense:
\end{oss}

\begin{defn} Let $\mathcal{F}$ an $\alpha-$sheaf over $X$. We call
\textit{support of} $\mathcal{F}$ the set
$$Supp\,\mathcal{F}:=\{x\in X\,|\,\mathcal{F}_{x}\neq
0\}.$$
\end{defn}

\begin{oss} If $\mathcal{F}$ is a coherent $\alpha-$sheaf, then
$Supp\,\mathcal{F}$ is a closed subset of $X$.
\end{oss}

\begin{oss}
\label{oss:x}Let $x\in X$ be a point. Then the skyscraper sheaf
$k(x)$ has a natural structure of $\alpha-$sheaf, for any
$\alpha\in Br'X$.
\end{oss}

We want also to recall some result on geometrical functors on
twisted sheaves. In particular we have:

\begin{prop}
\label{prop:funct} Let $X,Y$ be two varieties, $\alpha,\alpha'\in
Br'X$, $\beta\in Br'Y$ and $f:X\longrightarrow Y$ a morphism. Then
we can define the following
functors:$$\mathcal{H}om(.,.):Mod(X,\alpha)\times
Mod(X,\alpha')\longrightarrow Mod(X,\alpha'\alpha^{-1})$$
$$\otimes:Mod(X,\alpha)\times
Mod(X,\alpha')\longrightarrow Mod(X,\alpha\alpha')$$
$$f^{*}:Mod(Y,\beta)\longrightarrow Mod(X,f^{*}\beta)$$
$$f_{*}:Mod(X,f^{*}\beta)\longrightarrow Mod(Y,\beta).$$
\end{prop}

\proof See \cite{Cal}, Proposition 1.2.10.\endproof

\begin{prop}
\label{prop:llrf}Let $X$ be a variety. Then $\alpha\in BrX$ if and
only if there is an $\alpha-$sheaf $\mathcal{E}$ that is locally
free of finite rank. In this case, we have that the sheaf
$\mathcal{A}:=\mathcal{E}nd_{X,\alpha}(\mathcal{E})$ has the
structure of Azumaya algebra on $X$, and that there is an
equivalence
$$Mod(X,\alpha)\stackrel{\sim}\longrightarrow Mod(\mathcal{A}),\,\,\,\,\,\,\,\,
\mathcal{F}\mapsto\mathcal{F}\otimes\mathcal{E}^{\vee},$$where
$Mod(\mathcal{A})$ is the category of $\mathcal{O}_{X}-$modules
which have the structure of right $\mathcal{A}-$module.
\end{prop}

\proof See \cite{Cal}, Theorem 1.3.5. and Theorem 1.3.7.\endproof

If we pass to derived categories and functors, we can show the
following proposition:

\begin{prop}
\label{prop:derfunct} Let $X,Y$ be two varieties,
$\alpha,\alpha'\in BrX$, $\beta\in BrY$ and $f:X\longrightarrow Y$
a morphism. Then we can define the following
functors:$$\otimes^{L}:D^{-}(X,\alpha)\times
D^{-}(X,\alpha')\longrightarrow D^{-}(X,\alpha\alpha')$$
$$R\mathcal{H}om(.,.):D(X,\alpha)^{o}\times
D^{+}(X,\alpha')\longrightarrow D(X,\alpha'\alpha^{-1})$$
$$Lf^{*}:D^{-}(Y,\beta)\longrightarrow D^{-}(X,f^{*}\beta)$$
$$Rf_{*}:D(QCoh(X,f^{*}\beta))\longrightarrow
D(QCoh(Y,\beta)).$$If $f$ is a proper morphism, we have
$$Rf_{*}:D(X,f^{*}\beta)\longrightarrow D(Y,\beta).$$Moreover, if
$X$ and $Y$ are smooth of finite dimension and $f$ is a proper
morphism, we have$$\otimes^{L}:D^{b}(X,\alpha)\times
D^{b}(X,\alpha')\longrightarrow D^{b}(X,\alpha\alpha')$$
$$R\mathcal{H}om(.,.):D^{b}(X,\alpha)^{o}\times
D^{b}(X,\alpha')\longrightarrow D^{b}(X,\alpha'\alpha^{-1})$$
$$Lf^{*}:D^{b}(Y,\beta)\longrightarrow D^{b}(X,f^{*}\beta)$$
$$Rf_{*}:D^{b}(X,f^{*}\beta)\longrightarrow
D^{b}(Y,\beta).$$
\end{prop}

\proof See \cite{Cal}, Theorem 2.2.4 and Theorem 2.2.6.\endproof

For relations among these derived functors, see \cite{Cal},
Section 2.3.

\section[Reconstruction of a variety]{Reconstruction of a variety}

In this section we will show that any variety $X$ can be recovered
from the category $Coh(X,\alpha)$ of $\alpha-$sheaves, where
$\alpha\in BrX$. The idea is that we can give a ringed space
structure to the set $E_{X,\alpha}$ of irreducible Serre
subcategories of finite type of $Coh(X,\alpha)$, and that it is,
in fact, a scheme isomorphic to $X$.

In order to do so, we will introduce the notion of Serre
subcategory of an abelian category. Using the fact that skyscraper
sheaves of points of $X$ are twisted, we will show that
(irreducible) Serre subcategories of finite type of
$Coh(X,\alpha)$ are in bijective correspondence with the
(irreducible) closed subsets of $X$. This allows us to put a
topology on $E_{X,\alpha}$, which recovers the topology of $X$.
The problem will be to give a good definition of a structure sheaf
on $E_{X,\alpha}$, such that we can get an isomorphism between
$E_{X,\alpha}$ and $X$.

\subsection[Serre subcategories of an abelian category]{Serre
subcategories of an abelian category}

\begin{defn} Let $\mathcal{A}$ be an abelian category. A subcategory
$\mathcal{I}$ of $\mathcal{A}$ is called \textit{Serre
subcategory} if for every short exact sequence in $\mathcal{A}$
$$0\longrightarrow A\longrightarrow B\longrightarrow
C\longrightarrow 0$$we have $B\in\mathcal{I}$ if and only if
$A,C\in\mathcal{I}$.\\We say that $\mathcal{I}$ is a \textit{Serre
subcategory of finite type} if it is a Serre subcategory of
$\mathcal{A}$ generated by an element $A\in\mathcal{I}$, that is,
$\mathcal{I}$ is the smallest Serre subcategory of $\mathcal{A}$
that contains $A$. Such an $A$ will be called a \textit{generator}
for $\mathcal{I}$.\\We say that $\mathcal{I}$ is an
\textit{irreducible Serre subcategory} if it is not generated (as
Serre subcategory) by two proper Serre subcategories.
\end{defn}

\begin{esem}
\label{esem:supp} Let $Coh_{Z}(X,\alpha)$ be the full subcategory
of $Coh(X,\alpha)$ whose objects have support contained in the
closed set $Z$ of $X$. Then it is easy to show that it is a Serre
subcategory of $Coh(X,\alpha)$.
\end{esem}

\begin{defn}
If $\mathcal{I}$ is a subcategory of $\mathcal{A}$, the
\textit{quotient category} $\mathcal{A}/\mathcal{I}$ is the
category which has the same objects as $\mathcal{A}$ and morphisms
are defined in this way: if $A,B\in\mathcal{A}$ we have
$$Hom_{\mathcal{A}/\mathcal{I}}(A,B)=\lim_{\longrightarrow}Hom_{\mathcal{A}}(A',B')$$where
$i:A'\hookrightarrow A$ is a sub-object of $A$ such that
$coker(i)\in\mathcal{I}$ and $p:B\twoheadrightarrow B'$ is a
quotient of $B$ such that $ker(p)\in\mathcal{I}$.
\end{defn}

If $\mathcal{I}$ is a Serre subcategory of $\mathcal{A}$, then
$\mathcal{A}/\mathcal{I}$ is an abelian category. We have the
following lemmas, very easy to show:

\begin{lem}
\label{lem:equiv} Let $\mathcal{A},\mathcal{B}$ be abelian
categories and $F:\mathcal{A}\longrightarrow\mathcal{B}$ an exact
functor that admits a fully faithful right adjoint. Then $\ker F$
is a Serre subcategory of $\mathcal{A}$ and the induced functor
$F:\mathcal{A}/\ker F\longrightarrow\mathcal{B}$ is an
equivalence.
\end{lem}

\begin{lem}
\label{lem:plfid} Let $\mathcal{A}$ be an abelian category,
$\mathcal{A}'$ a full abelian subcategory of $\mathcal{A}$ and
$\mathcal{I}$ a Serre subcategory of $\mathcal{A}$. Suppose that
for every $M\in\mathcal{A}'$, $N\in\mathcal{I}$ a sub-object or a
quotient of $M$, we have that $N\in\mathcal{I}\cap\mathcal{A}'$.
Then the induced functor
$$i:\mathcal{A}'/\mathcal{I}\cap\mathcal{A}'\longrightarrow\mathcal{A}/\mathcal{I}$$is
fully faithful.
\end{lem}

Now, let $X$ be a variety over a field $k$, $\alpha\in BrX$ and
$Z$ a closed subset of $X$. Let $U=X\setminus Z$ and
$j_{U}:U\longrightarrow X$ the corresponding open immersion. We
have $$j_{U}^{*}:QCoh(X,\alpha)\longrightarrow
QCoh(U,\alpha_{|U})$$that is an exact functor with a fully
faithful right adjoint $j_{U*}$. Using Lemma \ref{lem:equiv} we
find that
$$j_{U}^{*}:QCoh(X,\alpha)/QCoh_{Z}(X,\alpha)\stackrel{\sim}\longrightarrow
QCoh(U,\alpha_{|U})$$is an equivalence and, using Lemma
\ref{lem:plfid}, that
$$j_{U}^{*}:Coh(X,\alpha)/Coh_{Z}(X,\alpha)\longrightarrow
Coh(U,\alpha_{|U})$$is fully faithful.

\begin{oss}
\label{oss:twist} If $\alpha=1$, so that $Coh(X,\alpha)\simeq
Coh(X)$, we can actually show that the functor $j_{U}^{*}$ above
is even an equivalence. Indeed, we know that every coherent sheaf
over an open subscheme $U$ of $X$ is restriction to $U$ of a
coherent sheaf on $X$. This is not clear in the case of twisted
sheaves: we will discuss this problem in Section 3.
\end{oss}

\subsection[Closed subsets and Serre subcategories]{Closed subsets and Serre
subcategories}

Let $X$ be a variety. In this section, we will show that for every
$\alpha\in BrX$, there is a bijective correspondence between
closed subsets of $X$ and Serre subcategories of finite type of
$Coh(X,\alpha)$. The main point here is the following:

\begin{prop}
\label{prop:gen} Let $X$ be a variety over a field $k$, $\alpha\in
BrX$ and $Z$ a closed subset of $X$. Then $Coh_{Z}(X,\alpha)$ is a
Serre subcategory of finite type of $Coh(X,\alpha)$. More
precisely, it is generated, as Serre subcategory, by any
$\alpha-$sheaf  $\mathcal{F}$ such that $Supp\,\mathcal{F}=Z$.
\end{prop}

\proof First, we have to show that such an $\alpha-$sheaf exists:
since $\alpha\in BrX$, there is a locally free $\alpha-$sheaf
$\mathcal{E}$ of finite rank over $X$ (see Proposition
\ref{prop:llrf}) so that its restriction to $Z$ (thought as an
$\alpha-$sheaf over $X$) has support equal to $Z$.  Now, choose
$\mathcal{F}\in Coh(X,\alpha)$ such that $Supp\,\mathcal{F}=Z$,
and write $\langle\mathcal{F}\rangle$ for the Serre subcategory of
$Coh(X,\alpha)$ generated by $\mathcal{F}$. It is clear that
$\langle\mathcal{F}\rangle\subseteq Coh_{Z}(X,\alpha)$, so that it
remains to show the opposite inclusion.

First, we reduce to $Z=X$: if we note $i:Z\longrightarrow X$ the
closed immersion of $Z$, let $\mathcal{I}$ be the Serre
subcategory of $Coh(Z,i^{*}\alpha)$ generated by
$i^{*}\mathcal{F}$. It is easy to show that
$i_{*}i^{*}\mathcal{F}\in\,\langle\mathcal{F}\rangle$, so that
$i_{*}\mathcal{I}\subseteq\langle\mathcal{F}\rangle$. If we have
that $\mathcal{I}=Coh(Z,i^{*}\alpha)$, we get
$Coh_{Z}(X,\alpha)=i_{*}\mathcal{I}\subseteq\,\langle\mathcal{F}\rangle$.
Using the same kind of arguments, we can even suppose $X$
irreducible.

We now proceed by induction over the dimension of $X$. The case of
dimension 0 is clear (here, every twist is trivial). Now, let us
suppose that $dimX=n$ and that the proposition is true for all
schemes of dimension smaller or equal to $n-1$. Let $Y$ be a
proper closed subscheme of $X$. By induction we have that
$Coh_{Y}(X,\alpha)\subseteq\,\langle\mathcal{F}\rangle$.

Now let $\mathcal{G}\in Coh(X,\alpha)$, $j:U\longrightarrow X$ an
open affine subscheme of $X$, and let $Y=X\setminus U$. If $U$ is
little enough, we can even suppose $j^{*}\mathcal{G}$ and
$j^{*}\mathcal{F}$ free of ranks $r$ and $s$ respectively (this is
possible, since $\alpha_{|U}\in BrU$). In this way we have an
isomorphism
$$\tilde{f}:j^{*}\mathcal{G}^{s}\stackrel{\sim}\longrightarrow
j^{*}\mathcal{F}^{r}$$in the category $Coh(U,\alpha_{|U})$. By
Lemma \ref{lem:plfid}, as we saw, the functor
$$j_{U}^{*}:Coh(X,\alpha)/Coh_{Y}(X,\alpha)\longrightarrow
Coh(U,\alpha_{|U})$$is fully faithful, so that the isomorphism
$\tilde{f}$ comes from an isomorphism
$f:\mathcal{G}^{s}\longrightarrow\mathcal{F}^{r}$ in
$Coh(X,\alpha)/Coh_{Y}(X,\alpha)$, that is $\ker f$ and co$\ker f$
are in $Coh_{Y}(X,\alpha)\subseteq\,\langle\mathcal{F}\rangle$.
But $\langle\mathcal{F}\rangle$ is a Serre subcategory, so
$\mathcal{G}\in\,\langle\mathcal{F}\rangle$.
\endproof

It is now easy to show the following:

\begin{cor}
\label{cor:bij} Let $X$ be a variety over a field $k$ and
$\alpha\in BrX$. There is a bijective correspondence between the
set $C$ of closed subsets of $X$ and the set $S$ of Serre
subcategories of finite type of $Coh(X,\alpha)$. In particular,
this induces a bijective correspondence between the points of $X$
and the set $E_{X,\alpha}$ of irreducible Serre subcategories of
finite type of $Coh(X,\alpha)$.
\end{cor}

\proof We can define$$i:C\longrightarrow
S,\,\,\,\,\,\,\,\,\,\,\,Z\mapsto Coh_{Z}(X,\alpha),$$and
$$j:S\longrightarrow
C,\,\,\,\,\,\,\,\,\,\,\,\mathcal{I}=\langle\mathcal{F}\rangle\mapsto
Supp\,\mathcal{F}.$$Now, $j$ is well defined: two generators of
the same Serre subcategory have the same support (this follows
from definition of generator and Remark \ref{oss:x}). In view of
Proposition \ref{prop:gen}, it is straightforward to show that
$i=j^{-1}$. Moreover, it is easy to show that $Z$ is irreducible
if and only if $Coh_{Z}(X,\alpha)$ is irreducible as Serre
subcategory of $Coh(X,\alpha)$. This gives the bijective
correspondence between the points of $X$ (which are the generic
points of irreducible closed sets of $X$) and the elements of
$E_{X,\alpha}$.
\endproof

\subsection[The reconstruction of a variety from $Coh(X,\alpha)$]
{The reconstruction of a variety from $Coh(X,\alpha)$}

We are now able to describe how one can recover the variety $X$
from $Coh(X,\alpha)$. As we saw in Corollary \ref{cor:bij}, we can
think the points of $X$ as the irreducible Serre subcategories of
finite type of $Coh(X,\alpha)$.

Let $E=E_{X,\alpha}$ be the set of irreducible Serre subcategories
of finite type of $Coh(X,\alpha)$. On $E$ we can define the
following topology: let $\mathcal{I}$ be a Serre subcategory of
finite type of $Coh(X,\alpha)$, and write
$$D(\mathcal{I}):=\{\mathcal{J}\in
E\,|\,\mathcal{J}\nsubseteq\mathcal{I}\}.$$It is easy to verify
that this family of subsets forms a topology over $E$ and that the
following morphism:
$$f:=f_{X,\alpha}:E\longrightarrow
X,\,\,\,\,\,\,\,\,\,f(\mathcal{J}=Coh_{\overline{\{x\}}}(X,\alpha))=x$$is
a homeomorphism (use Proposition \ref{prop:gen} and Corollary
\ref{cor:bij}). More precisely, if $Z$ is a closed subset of $X$,
$U=X\setminus Z$ and $\mathcal{I}=Coh_{Z}(X,\alpha)$, then $f$
gives a bijective correspondence between $D(\mathcal{I})$ and $U$.
In this way, we have shown that we can recover the topological
space underlying
$X$ from $Coh(X,\alpha)$.\\

\indent It remains to define a structure sheaf on $E$, in order to
make $f$ an isomorphism of schemes. Let us recall the notion of
center of a category.

\begin{defn} Let $\mathcal{A}$ be a category. We call \textit{center
of} $\mathcal{A}$ the ring $Z(\mathcal{A})$ of endomorphisms of
the identity functor of $\mathcal{A}$:
$Z(\mathcal{A}):=End_{\mathcal{A}}(id_{\mathcal{A}})$.
\end{defn}

It is very easy to show the following lemma:

\begin{lem}
\label{lem:z} Let $A$ be a ring with unity, $Z(A)$ his center, and
$Mod_{ft}(A)$ the category of modules of finite type over $A$. The
canonical morphism
$$Z(A)\longrightarrow Z(Mod_{ft}(A)),\,\,\,\,\,\,a\mapsto\cdot a$$is
an isomorphism of commutative rings.
\end{lem}

Using the center of a category and the notations we used above, we
can define the following sheaf on $E$:
$$\mathcal{O}_{E}(D(\mathcal{I}))=Z(Coh(U,\alpha_{|U})),$$and the
morphism of sheaves $f^{\natural}:\mathcal{O}_{X}\longrightarrow
f_{*}\mathcal{O}_{E}$ which is given over every open set $U$ of
$X$ by
$$f^{\natural}(U):\mathcal{O}_{X}(U)\longrightarrow
Z(Coh(U,\alpha_{|U})),\,\,\,\,\,\,\,\,\,\,s\mapsto\cdot s.$$In
this way we have given to $E$ the structure of ringed space. We
have now the following theorem, which shows the first part of
Theorem \ref{thm:gab}:

\begin{thm}
\label{thm:gabriel} The morphism
$(f,f^{\natural}):E\longrightarrow X$ is an isomorphism of ringed
spaces over $k$. In particular, $E$ is a $k-$variety which depends
only on $Coh(X,\alpha)$.
\end{thm}

\proof We only need to show that $f^{\natural}$ is an isomorphism
of rings. It suffices to show that it is an isomorphism on open
affine subschemes of $X$. So, let's take $U=Spec\,A$ an open
affine subscheme of $X$.

As we want to show that $f^{\natural}(U)$ is an isomorphism, we
begin by studying the ring $Z(Coh(U,\alpha_{|U}))$. Since
$\alpha_{|U}\in BrU$, from Proposition \ref{prop:llrf} we know
that there is a locally free $\alpha-$sheaf $\mathcal{E}$ of rank
$r$, and that there is an equivalence of categories
$$Mod(U,\alpha_{|U})\stackrel{\sim}\longrightarrow
Mod(\mathcal{E}nd_{U,\,\alpha_{|U}}(\mathcal{E})),\,\,\,\,\,\,
\mathcal{F}\mapsto\mathcal{F}\otimes\mathcal{E}^{\vee}.$$If we
note $\mathcal{A}=\mathcal{E}nd_{U,\,\alpha_{|U}}(\mathcal{E})$,
this is an Azumaya algebra on $X$ (so that, in particular, it is
an $\mathcal{O}_{X}-$module). Looking at this equivalence, we can
see that it sends any (quasi) coherent $\alpha_{|U}-$sheaf to a
(quasi) coherent sheaf which has the structure of right
$\mathcal{A}-$module, so that we get an equivalence between
$Coh(U,\alpha_{|U})$ and the full subcategory $\mathcal{C}$ of
$Mod(\mathcal{A})$ whose objects are coherent sheaves with the
structure of right $\mathcal{A}-$module. Since we are on an affine
scheme, taking global sections we get an equivalence between
$\mathcal{C}$ and $Mod_{ft}(End_{U,\,\alpha_{|U}}(\mathcal{E}))$,
so that we finally have the isomorphisms
$$Z(Coh(U,\alpha_{|U}))\simeq
Z(Mod_{ft}(End_{U,\alpha_{|U}}(\mathcal{E})))\simeq
Z(End_{U,\alpha_{|U}}(\mathcal{E}))$$where the second isomorphism
follows from Lemma \ref{lem:z}.

We are now reduced to study the center of the ring of
endomorphisms of $\mathcal{E}$ as an $\alpha_{|U}-$sheaf. Let's
see what this ring looks like. As $U=Spec\,A$, we can find
$f_{1},...,f_{n}\in A$ such that $U=\bigcup_{i=1}^{n}D(f_{i})$. We
represent $\mathcal{E}$ over this open covering, so that
$\mathcal{E}=(\mathcal{E}_{i},\varphi_{ij})$, where
$\mathcal{E}_{i}$ is an locally free $A_{f_{i}}-$module of rank
$r$ for every $i$, and
$$\varphi_{ij}:\mathcal{E}_{j|D(f_{i}f_{j})}\stackrel{\sim}
\longrightarrow\mathcal{E}_{i|D(f_{i}f_{j})}$$is an isomorphism of
$A_{f_{i}f_{j}}-$modules for every $i,j$, verifying the conditions
of Definition \ref{defn:twisted}. We can even take $D(f_{i})$
small enough such that, for every $i$, $\mathcal{E}_{i}\simeq
(A_{f_{i}})^{r}$. With this choice, we see that $\varphi_{ij}$ can
be thought as an isomorphism of $(A_{f_{i}f_{j}})^{r}$, that is a
matrix $B_{ij}\in GL_{r}(A_{f_{i}f_{j}})$.

Now, to give an endomorphism $M$ of $\mathcal{E}$ is to give a
family $(M_{1},...,M_{n})$, where $M_{i}\in
End_{A_{f_{i}}}((A_{f_{i}})^{r})=M_{r}(A_{f_{i}})$ is a square
matrix of rank $r$ with elements in $A_{f_{i}}$, such that for
every $i,j$ we have $B_{ij}M_{i}=M_{j}B_{ij}$ in the ring
$M_{r}(A_{f_{i}f_{j}})$. In conclusion, we have
$$End_{U,\,\alpha_{|U}}(\mathcal{E})=\{(M_{1},...,M_{n})\,|\,M_{i}\in
M_{r}(A_{f_{i}}),\, B_{ij}M_{i}=M_{j}B_{ij}\in
M_{r}(A_{f_{i}f_{j}})\}$$where sum and multiplication are the
obvious ones. It is easy to see that the morphism
$f^{\natural}(U)$ is now given by $$A\longrightarrow
End_{U,\,\alpha_{|U}}(\mathcal{E}),\,\,\,\,\,\,a\mapsto(diag(a),...,diag(a)).$$It
is quite clear that this map is injective: if $a,b\in A$ have the
same image, this means that $a=b$ in $A_{f_{i}}$ for every $i$.
Since $a,b$ are global section of the sheaf associated to $A$,
this means that $a=b$ in $A$.

It remains to show the surjectivity. Let $M=(M_{1},...,M_{n})$ be
an element in $Z(End_{U,\,\alpha_{|U}}(\mathcal{E}))$. This means
that for every $(N_{1},...,N_{n})\in
End_{U,\,\alpha_{|U}}(\mathcal{E})$ we have
$$(M_{1}N_{1},...,M_{n}N_{n})=(N_{1}M_{1},...,N_{n}M_{n}),$$that
is, for every $i$, $M_{i}\in Z(M_{r}(A_{f_{i}}))\simeq A_{f_{i}}$,
so that there is $b_{i}\in A_{f_{i}}$ such that
$M_{i}=diag(b_{i})$. In particular, $M_{i}\in
Z(M_{r}(A_{f_{i}f_{j}}))$ for every $i,j$. The condition
$B_{ij}M_{i}=M_{j}B_{ij}$ in the ring $M_{r}(A_{f_{i}f_{j}})$
gives $diag(b_{i}-b_{j})B_{ij}=0$. Since $B_{ij}$ is invertible,
this tells that $b_{i}=b_{j}$ in $A_{f_{i}f_{j}}$ for every $i,j$.
Since $b_{i}$ is a section on $D(f_{i})$ of the sheaf associated
to $A$, we get an element $a\in A$ such that $b_{i}=a$ in
$A_{f_{i}}$ for every $i$. This tells us that
$(M_{1},...,M_{n})=(diag(a),...,diag(a))$ for a (unique) $a\in
A$.\endproof

\section[Isomorphism induced by an equivalence]{Isomorphism
induced by an equivalence}

We have shown a generalization of the first part of Gabriel's
theorem to twisted coherent sheaves, as we wanted at the
beginning, so, the next question is if any equivalence between
$Coh(X,\alpha)$ and $Coh(Y,\beta)$ gives rise to an isomorphism
between $X$ and $Y$.

First of all, we can show that this problem can be reduced to the
following one: let $X$ be a variety, $\alpha\in BrX$, $Z$ a closed
subset, $U=X\setminus Z$ and $j_{U}:U\longrightarrow X$ the open
immersion. Is the functor
$$j_{U}^{*}:Coh(X,\alpha)/Coh_{Z}(X,\alpha)\longrightarrow
Coh(U,\alpha_{|U})$$an equivalence?

If we have a positive answer for every open subscheme $U$ we'll
say that $(X,\alpha)$ \textit{satisfies the restriction
condition}.

So, let $(Y,\beta)$ be another variety which satisfies the
restriction condition,
and$$F:Coh(X,\alpha)\stackrel{\sim}\longrightarrow
Coh(Y,\beta)$$an equivalence. It is trivial to show that if
$\mathcal{I}$ is a(n irreducible) Serre subcategory of finite type
of $Coh(X,\alpha)$, then $F(\mathcal{I})$ is a(n irreducible)
Serre subcategory of finite type of $Coh(Y,\beta)$. This gives a
bijective correspondence
$$f_{F}:X\longrightarrow Y,\,\,\,\,\,\,\,\,f_{F}(x)=
f_{Y,\beta}(F(f_{X,\alpha}^{-1}(x))),$$where $f_{X,\alpha}$ (resp.
$f_{Y,\beta}$) is the isomorphism between $E_{X,\alpha}$ and $X$
(resp. between $E_{Y,\beta}$ and $Y$) we defined in the previous
section. It is also easy to show that $f_{F}$ is an homeomorphism
and even an isomorphism of schemes: $U$ is an open subscheme of
$X$, we have that $f_{F}$ induces a bijective correspondence
between $U$ and $W:=f_{Y,\beta}(D(F(f_{X,\alpha}^{-1}(U))))$.
Moreover, since $F$ is an equivalence, we have that
$$F:Coh(X,\alpha)/Coh_{X\setminus U}(X,\alpha)\stackrel{\sim}\longrightarrow
Coh(Y,\beta)/Coh_{Y\setminus W}(Y,\beta)$$is an equivalence (this
follows easily from Lemma \ref{lem:equiv}). Using this and the
fact that $(X,\alpha)$ and $(Y,\beta)$ verify the restriction
condition, it's trivial to show that
$$j_{W}^{*}\circ F\circ (j_{U}^{*})^{-1}:Z(Coh(U,\alpha_{|U}))
\stackrel{\sim}\longrightarrow Z(Coh(W,\beta_{|W}))$$is an
isomorphism, that is, we get an isomorphism
$g:\mathcal{O}_{E_{X,\alpha}}\stackrel{\sim}\longrightarrow\mathcal{O}_{E_{Y,\beta}}$.
Now, from Theorem \ref{thm:gabriel}, it's obvious that
$$f_{F}^{\natural}:=f_{X,\alpha}^{\natural}\circ g^{-1}\circ f_{Y,\beta}^{\natural\,-1}
:\mathcal{O}_{Y}\stackrel{\sim}\longrightarrow
f_{F*}\mathcal{O}_{X}$$is an isomorphism. In conclusion, we have
shown the following:

\begin{thm}
\label{thm:if} Let $X,Y$ be two varieties over a field $k$,
$\alpha\in BrX$ and $\beta\in BrY$ which verify the restriction
condition above. Then any equivalence
$$F:Coh(X,\alpha)\stackrel{\sim}\longrightarrow Coh(Y,\beta)$$induces an
isomorphism of varieties $f:X\stackrel{\sim}\longrightarrow Y.$
\end{thm}

We are now reduced to study when a couple $(X,\alpha)$ verifies
the restriction condition. Actually, we can show it only when $X$
is a smooth variety. To approach the problem, we get into the
domain of derived category, where we can use perfect and compact
objects in $D(QCoh(X,\alpha))$.

\subsection[Thick subcategories and compact objects of a triangulated
category]{Thick subcategories and compact objects of a
triangulated category}

In this section we introduce the notions of thick subcategory and
Bousfield subcategory of a triangulated category. The main
references here will be \cite{Ro1} and \cite{Ro2}. Let
$\mathcal{T}$ be a triangulated category.

\begin{defn} We say that a subcategory $\mathcal{I}$ of
$\mathcal{T}$ is \textit{thick} (or \textit{épaisse}) if it is a
triangulated subcategory such that for every $M,N\in\mathcal{T}$,
if $M\oplus N\in\mathcal{I}$ then $M,N\in\mathcal{I}$.
\end{defn}

\noindent If $\mathcal{I}$ is a thick subcategory of
$\mathcal{T}$, we have that the quotient category
$\mathcal{T}/\mathcal{I}$ is again triangulated. It is clear that
we have a(n essentially surjective) functor
$$j^{*}:\mathcal{T}\longrightarrow\mathcal{T}/\mathcal{I}.$$

\begin{defn} A thick subcategory $\mathcal{I}$ of $\mathcal{T}$ is
called \textit{Bousfield subcategory} if $j^{*}$ admits a right
adjoint, which will be noted $j_{*}$.
\end{defn}

Now, let $\mathcal{I}$ be a full triangulated subcategory of
$\mathcal{T}$. We define the following subcategories of
$\mathcal{T}$:
\begin{enumerate}
\item If $\mathcal{T}$ admits infinite direct sums,
$\overline{\mathcal{I}}$ will be smallest thick subcategory of
$\mathcal{T}$ which contains $\mathcal{I}$ and which is stable for
infinite direct sums;\item $\langle\mathcal{I}\rangle$ is the
smallest thick subcategory of $\mathcal{T}$ which contains
$\mathcal{I}$;\item $\mathcal{I}^{\bot}$ is the subcategory of
objects $C\in\mathcal{T}$ such that for every
$D\in\,\langle\mathcal{I}\rangle$ we have
$Hom_{\mathcal{T}}(D,C)=0$;\item if
$\mathcal{I}_{1},\mathcal{I}_{2}$ are two full subcategories of
$\mathcal{T}$ we define $\mathcal{I}_{1}\ast\mathcal{I}_{2}$ as
the subcategory of $\mathcal{T}$ whose objects $M$ are such that
there is a distinguished triangle $$M_{1}\longrightarrow
M\longrightarrow M_{2}\rightsquigarrow$$where
$M_{i}\in\mathcal{I}_{i}$, and
$\mathcal{I}_{1}\diamond\mathcal{I}_{2}=\langle\mathcal{I}_{1}
\ast\mathcal{I}_{2}\rangle$; \item
$\langle\mathcal{I}\rangle_{0}=0$, and by induction over $i$ we
define
$\langle\mathcal{I}\rangle_{i}=\langle\mathcal{I}\rangle_{i-1}
\diamond\langle\mathcal{I}\rangle$ and
$\langle\mathcal{I}\rangle_{\infty}=\bigcup_{i\geq
0}\langle\mathcal{I}\rangle_{i}$. In particular, if $\mathcal{T}$
admits infinite direct sums, we have
$\langle\mathcal{I}\rangle_{\infty}=\overline{\langle\mathcal{I}\rangle}$.
\end{enumerate}

\begin{lem}
If $\mathcal{I}_{1}$ et $\mathcal{I}_{2}$ are two Bousfield
subcategories of $\mathcal{T}$, then
$\mathcal{I}_{1}\cap\mathcal{I}_{2}$ and
$\langle\mathcal{I}_{1}\cup\mathcal{I}_{2}\rangle_{\infty}$ are
Bousfield subcategories of $\mathcal{T}$.
\end{lem}

\proof See \cite{Ro2}, Lemma 5.8.\endproof

\begin{defn} Un object $C\in\mathcal{T}$ is said to be
\textit{compact} in $\mathcal{T}$ if for every family
$\mathcal{E}$ of objects of $\mathcal{T}$, the canonical morphism
$$\bigoplus_{E\in\mathcal{E}}Hom(C,E)\longrightarrow
Hom(C,\bigoplus_{E\in\mathcal{E}}E)$$is an isomorphism. We will
note $\mathcal{T}^{c}$ the full subcategory of $\mathcal{T}$ of
compact objects.
\end{defn}

We can now recall some result which will be used later.

\begin{lem}
\label{lem:equivtri} Let $\mathcal{T},\mathcal{T}'$ be two
triangulated categories and
$F:\mathcal{T}\longrightarrow\mathcal{T}'$ an exact functor which
admits a fully faithful right adjoint. Then $\ker F$ is a thick
subcategory of $\mathcal{T}$ and the induced functor
$F:\mathcal{T}/\ker F\longrightarrow\mathcal{T}'$ is an
equivalence.
\end{lem}

\begin{lem}
\label{lem:plfidtri} Let $\mathcal{T}$ be a triangulated category,
$\mathcal{I}$ a thick subcategory of $\mathcal{T}$ and
$\mathcal{T}'$ a full triangulated subcategory of $\mathcal{T}$.
If for every $C\in\mathcal{T}'$, $D\in\mathcal{I}$, every morphism
from $C$ to $D$ factorizes by an object of
$\mathcal{I}\cap\mathcal{T}'$, we have that the induced functor
$$i:\mathcal{T}'/\mathcal{I}\cap\mathcal{T}'\longrightarrow\mathcal{T}/\mathcal{I}$$is
fully faithful.
\end{lem}

\begin{lem}
\label{lem:neeman} Let $\mathcal{T}$ be a triangulated category
that admits infinite direct sums and $\mathcal{I}$ a thick
subcategory of $\mathcal{T}^{c}$. Then any morphism from a compact
object of $\mathcal{T}$ to an object of $\overline{\mathcal{I}}$
factorizes by an element of $\mathcal{I}$. In particular
$\mathcal{T}^{c}\cap\overline{\mathcal{I}}=\mathcal{I}$. Moreover,
$\overline{\mathcal{I}}=\mathcal{T}$ if and only if
$\mathcal{I}^{\bot}=0$.
\end{lem}

\proof This is proved in \cite{BN}.\endproof

\begin{lem}
\label{lem:mayer} Let $\mathcal{I}_{1},\mathcal{I}_{2}$ be two
Bousfield subcategories of $\mathcal{T}$ with
$\mathcal{I}_{1}\cap\mathcal{I}_{2}=0$. Then for every object
$D\in\mathcal{T}$ there is a distinguished triangle:
$$D\longrightarrow j_{1*}j_{1}^{*}D\oplus
j_{2*}j_{2}^{*}D\longrightarrow j_{\cup
*}j_{\cup}^{*}D\longrightarrow D[1]$$which is called
\textit{Mayer-Vietoris triangle for} $D$, where
$j_{i}^{*}:\mathcal{T}\longrightarrow\mathcal{T}/\mathcal{I}_{i}$
and $j_{\cup}^{*}:\mathcal{T}\longrightarrow\mathcal{T}/\langle
\mathcal{I}_{1}\cup\mathcal{I}_{2}\rangle_{\infty}$ are the
projection functors.
\end{lem}

\proof See \cite{Ro2}, Proposition 5.10.\endproof

\subsection[Perfect objects in $D(QCoh(X,\alpha))$]{Perfect objects in
$D(QCoh(X,\alpha))$}

Let $X$ be a variety over a field $k$ and $\alpha\in BrX$. If $Z$
is a closed subset of $X$, it is easy to see that
$D_{Z}(QCoh(X,\alpha))$ is a thick subcategory of
$D(QCoh(X,\alpha))$. Moreover, using Lemma \ref{lem:equivtri} we
can see that if $U=X\setminus Z$ and $j_{U}$ is the open immersion
of $U$ in $X$, then
$$j_{U}^{*}:D(QCoh(X,\alpha))/D_{Z}(QCoh(X,\alpha))\stackrel{\sim}\longrightarrow
D(QCoh(U,\alpha_{|U}))$$is an equivalence, since we have the
functor $Rj_{U*}$ which is fully faithful and is right adjoint to
$j_{U}^{*}$. In this way we have also shown that
$D_{Z}(QCoh(X,\alpha))$ is a Bousfield subcategory of
$D(QCoh(X,\alpha))$.

Now, let $Z_{1},Z_{2}$ be two closed subsets of $X$ such that
$Z_{1}\cap Z_{2}=\emptyset$, and let $U_{i}=X\setminus Z_{i}$,
$U_{12}=U_{1}\cap U_{2}$. It's clear that
$$D_{Z_{1}}(QCoh(X,\alpha))\cap D_{Z_{2}}(QCoh(X,\alpha))=0,$$and
it's easy to show that
$$\langle D_{Z_{1}}(QCoh(X,\alpha))\cup
D_{Z_{2}}(QCoh(X,\alpha))\rangle_{\infty}=D_{Z_{1}\cup
Z_{2}}(QCoh(X,\alpha)).$$Using Lemma \ref{lem:mayer}, if $D\in
D(QCoh(X,\alpha))$, there is a distinguished triangle
\begin{equation}
\label{eq:triangle}D\longrightarrow Rj_{1*}j_{1}^{*}D\oplus
Rj_{2*}j_{2}^{*}D\longrightarrow
Rj_{12*}j_{12}^{*}D\longrightarrow D[1]
\end{equation}
where $j_{i}:U_{i}\longrightarrow X$ and
$j_{12}:U_{12}\longrightarrow X$.

We give the following definition:

\begin{defn} An object $C\in D(QCoh(X,\alpha))$ is called
\textit{perfect} if it is locally quasi-isomorphic to a bounded
complex of locally free $\alpha-$sheaves of finite rank. We will
denote $Perf(X,\alpha)$ the subcategory of perfect objects in
$D(QCoh(X,\alpha))$.
\end{defn}

$Perf(X,\alpha)$ is a (non empty) thick subcategory of
$D^{b}(X,\alpha)$. Moreover, if $X$ is a smooth variety, we have
$Perf(X,\alpha)=D^{b}(X,\alpha)$ (see \cite{Cal}, Lemma 2.1.4 and
Proposition 2.1.8).

We have the following theorem, which will be basic for what will
follow.

\begin{thm}
\label{thm:perf} Let $X$ be a variety over a field $k$, $\alpha\in
BrX$. Then we have $Perf(X,\alpha)=D(QCoh(X,\alpha))^{c}$.
\end{thm}

\proof The proof will proceed by induction on the minimal number
of open affine subschemes which cover $X$. We start with
$X=Spec\,A$. Since we know that $\alpha\in BrX$, there is a
locally free $\alpha-$sheaf $\mathcal{E}$ of finite rank over $X$.

First of all, $\mathcal{E}$ is compact: let
$\{\mathcal{F}_{i}\}_{i\in I}$ a family of complexes of
quasi-coherent $\alpha-$sheaves, and let
$\mathcal{E}^{\vee}=R\mathcal{H}om(\mathcal{E},\mathcal{O}_{X})$
be the dual complex of $\mathcal{E}$. We have
$$\bigoplus_{i\in I}Hom_{X,\alpha}(\mathcal{E},\mathcal{F}_{i})=
\bigoplus_{i\in
I}Hom_{X}(\mathcal{O}_{X},\mathcal{E}^{\vee}\otimes\mathcal{F}_{i}),$$and
since $\mathcal{O}_{X}$ is compact in $D(QCoh(X))$ (this is
trivial), we have
$$\bigoplus_{i\in I}Hom_{X}(\mathcal{O}_{X},\mathcal{E}^{\vee}\otimes\mathcal{F}_{i})=
Hom_{X}(\mathcal{O}_{X},\bigoplus_{i\in
I}\mathcal{E}^{\vee}\otimes\mathcal{F}_{i})=
Hom_{X,\alpha}(\mathcal{E}^{\vee},\bigoplus_{i\in
I}\mathcal{F}_{i}),$$so that $Perf(X,\alpha)\subseteq
D(QCoh(X,\alpha))^{c}$. We have even that $Perf(X,\alpha)$ is a
thick subcategory of $D(QCoh(X,\alpha))^{c}$.

Now we can show that
$\overline{Perf(X,\alpha)}=D(QCoh(X,\alpha))$: using Lemma
\ref{lem:neeman}, it suffices to show that
$Perf(X,\alpha)^{\bot}=0$. So, let $C\in Perf(X,\alpha)^{\bot}$.
Since $\mathcal{E}$ is perfect, we have that
$$0=RHom_{X,\alpha}(\mathcal{E},C)=RHom_{X}(\mathcal{O}_{X},\mathcal{E}^{\vee}\otimes C)$$that
is $\mathcal{H}^{i}(\mathcal{E}^{\vee}\otimes C)=0$ for every $i$.
This implies clearly $\mathcal{E}^{\vee}\otimes C=0$, that is
$C=0$. Using again Lemma \ref{lem:neeman}, we have
$Perf(X,\alpha)=D(QCoh(X,\alpha))^{c}$.

Now suppose that $X=U_{1}\cup U_{2}$, where $U_{1}$ is affine and
$U_{2}$ verifies the theorem. Let $C,D\in D(QCoh(X,\alpha))$.
Using the Mayer-Vietoris triangle (\ref{eq:triangle}), it is easy
to show that $C$ is compact if and only if $j_{1}^{*}C$,
$j_{2}^{*}C$ and $j_{12}^{*}C$ are.

Since $U_{1}$, $U_{2}$ and $U_{12}$ verify the theorem by
induction, we have that $j_{1}^{*}C$, $j_{2}^{*}C$ and
$j_{12}^{*}C$ are compact if and only if they are perfect, that
is, there are $E_{i}\in D(QCoh(U_{i},\alpha_{|U_{i}}))$ and
$E_{12}\in D(QCoh(U_{12},\alpha_{|U_{12}}))$ bounded complexes of
locally free twisted sheaves of finite rank such that $j_{i}^{*}C$
is locally quasi-isomorphic to $E_{i}$ and $j_{12}^{*}C$ is
locally quasi-isomorphic to $E_{12}$. This means that
$E_{1|U_{12}}\simeq E_{12}\simeq E_{2|U_{12}}$, that is, we can
glue $E_{1}$ and $E_{2}$ over $E_{12}$, obtaining a locally free
$\alpha-$sheaf of finite rank on $X$, locally quasi-isomorphic to
$C$.
\endproof

Let us denote $Perf_{Z}(X,\alpha)=Perf(X,\alpha)\cap
D_{Z}(QCoh(X,\alpha))$.

\begin{defn} We define the group $K_{0}(\mathcal{T})$ of a triangulated
category $\mathcal{T}$ as the quotient of the free abelian group
generated by the objects in $\mathcal{T}$ by the relation
$[N]=[M]+[L]$ if there is a distinguished triangle
$$M\longrightarrow N\longrightarrow L\longrightarrow M[1].$$In
particular, we will denote $K_{0}(X,\alpha)=K_{0}(Perf(X,\alpha))$
and, if $Z$ is a closed subset of $X$,
$K_{0,Z}(X,\alpha)=K_{0}(Perf_{Z}(X,\alpha))$.
\end{defn}

We have the following lemma, due to Thomason.

\begin{lem}
\label{lem:k0} Let $\mathcal{T}$ be a triangulated category. The
correspondence which sends a full triangulated subcategory
$\mathcal{I}$ of $\mathcal{T}$ such that
$\overline{\mathcal{I}}=\mathcal{T}$ to the image of
$K_{0}(\mathcal{I})$ in $K_{0}(\mathcal{T})$ of the group morphism
induced by the inclusion $i:\mathcal{I}\longrightarrow\mathcal{T}$
is bijective.
\end{lem}

Thanks to this lemma, and the others we stated in section 3.1, we
can show the following:

\begin{thm}
\label{thm:obstruction} Let $X$ be a variety over a field $k$,
$\alpha\in BrX$. Let $Y,Z$ two closed subsets of $X$,
$U=X\setminus Z$ and $j_{U}$ the corresponding open immersion.
Then the functor
$$j_{U}^{*}:Perf_{Y}(X,\alpha)/Perf_{Z\cap
Y}(X,\alpha)\longrightarrow Perf_{U\cap Y}(U,\alpha_{|U})$$is
fully faithful. Moreover, an object $\mathcal{F}\in Perf_{U\cap
Y}(X,\alpha_{|U})$ is restriction to $U$ of an object in
$Perf_{Y}(X,\alpha)$ if and only if $[\mathcal{F}]\in K_{0,U\cap
Y}(U,\alpha_{|U})$ is restriction of a class in
$K_{0,Y}(X,\alpha)$.
\end{thm}

\proof Suppose that for every variety $X$, $\alpha\in BrX$ and $Y$
closed subset of $X$, we have that
$$\overline{Perf_{Y}(X,\alpha)}=D_{Y}(QCoh(X,\alpha)).$$Using this
condition, Lemma \ref{lem:plfidtri} and Lemma \ref{lem:neeman} we
can easily see that the functor $j_{U}^{*}$ is fully faithful.
Now, let $\mathcal{I}$ be the essential image of $j_{U}^{*}$. This
is a full triangulated subcategory of $Perf_{U\cap
Y}(U,\alpha_{|U})$. Since we have supposed that
$\overline{Perf_{Y}(X,\alpha)}=D_{Y}(QCoh(X,\alpha))$, we have
$\overline{\mathcal{I}}=D_{U\cap Y}(QCoh(U,\alpha_{|U}))$, and
from Lemma \ref{lem:neeman} we have that $\mathcal{I}$ generates
$Perf_{U\cap Y}(X,\alpha)$ as thick subcategory. Now the theorem
follows from Lemma \ref{lem:k0}.

We have now to show that
$\overline{Perf_{Y}(X,\alpha)}=D_{Y}(QCoh(X,\alpha))$ for every
variety $X$, $\alpha\in BrX$ and $Y$ a closed subset of $X$. Using
Lemma \ref{lem:neeman}, it suffices to show that
$Perf_{Y}(X,\alpha)^{\bot}=0$. We will proceed by induction on the
minimal number of open affine subschemes that cover $X$. So, let
$C\in Perf_{Y}(X,\alpha)^{\bot}$.

Let $X=Spec\,A$, so that $Y$ will be the closed subset
corresponding to the ideal of $A$ generated by $r$ elements
$f_{1},...,f_{r}$, or $Y=\emptyset$, that is $r=0$. Let
$\mathcal{E}$ a locally free $\alpha-$sheaf on $X$ and
$$G_{r}=\bigotimes_{i=1}^{r}(0\longrightarrow\mathcal{O}_{X}\otimes
\mathcal{E}^{\vee}\stackrel{\cdot f_{i}\otimes
id}\longrightarrow\mathcal{O}_{X}\otimes
\mathcal{E}^{\vee}\longrightarrow 0)\in Perf_{Y}(X,\alpha^{-1}).$$
It is easy to show that $C=0$ if and only if $G_{r}\otimes C=0$
(here we have that $G_{r}\otimes C$ is a sheaf, so we can use the
same argument in the proof of Lemme 2.10 in \cite{Ro1}). Now,
since $G_{r}\in Perf_{Y}(X,\alpha^{-1})$, we have that
$G_{r}^{\vee}\in Perf_{Y}(X,\alpha)$, so that
$$0=RHom_{X,\alpha}(G_{r}^{\vee},C)=RHom_{X}(\mathcal{O}_{X},G_{r}\otimes C)$$that
is $\mathcal{H}^{i}(G_{r}\otimes C)=0$ for every $i$, and so
$G_{r}\otimes C=0$, which implies $C=0$.

Now, let $X=U_{1}\cup U_{2}$, where $U_{1}$ is affine and $U_{2}$
verifies the theorem. Moreover let $Z_{i}=X\setminus U_{i}$. Let
$C\in Perf_{Y}(X,\alpha)^{\bot}$, $j_{i}$ the open immersion of
$U_{i}$ in $X$, $j_{12}$ the open immersion of $U_{12}=U_{1}\cap
U_{2}$ in $X$.

First we show that the adjonction morphism
$\gamma:C\longrightarrow Rj_{2*}j_{2}^{*}C$ is a
quasi-isomorphism. Let $D\in Perf_{Y\cap
Z_{2}}(U_{1},\alpha_{|U_{1}})$. Since $Y\cap Z_{2}\subseteq
U_{1}$, the functor $Rj_{1*}$ induces the two following
equivalences:
\begin{equation}
\label{eq:equiv} D_{Y\cap
Z_{2}}(QCoh(U_{1},\alpha_{|U_{1}}))\stackrel{Rj_{1*}}\longrightarrow
D_{Y\cap Z_{2}}(QCoh(X,\alpha))
\end{equation}
and
\begin{equation}
\label{eq:equiv2} Perf_{Y\cap
Z_{2}}(U_{1},\alpha_{|U_{1}})\stackrel{Rj_{1*}}\longrightarrow
Perf_{Y\cap Z_{2}}(X,\alpha).
\end{equation}
Using this, we have that $Rj_{1*}D\in Perf_{Y\cap
Z_{2}}(X,\alpha)\subseteq Perf_{Y}(X,\alpha)$, so that
$Hom(Rj_{1*}D,C)=0$. Moreover
$$Hom(Rj_{1*}D,Rj_{2*}j_{2}^{*}C)=Hom(j_{2}^{*}Rj_{1*}D,j_{2}^{*}C)=0.$$
If $C'$ is the cocone of $\gamma$, applying the functor
$Hom(Rj_{1*}D,.)$ to the distinguished triangle
$$C'\longrightarrow C\stackrel{\gamma}\longrightarrow
Rj_{2*}j_{2}^{*}C\longrightarrow C'[1]$$we get
$Hom(Rj_{1*}D,C')=0$. Now, we know that $C'\in D_{Y\cap
Z_{2}}(QCoh(X,\alpha))$. Using equivalence (\ref{eq:equiv}), we
get $C''\in D_{Y\cap Z_{2}}(QCoh(U_{1},\alpha_{|U_{1}}))$ with the
property that $Rj_{1*}C''\simeq C'$. In this way we have that
$Hom(D,C'')=0$ for every $D\in Perf_{Y\cap
Z_{2}}(U_{1},\alpha_{|U_{1}})$. But $U_{1}$ is affine, so the
first part of the proof says $C''=0$, and then $C'=0$, that is,
$\gamma$ is a quasi-isomorphism.

Now let $E'\in Perf_{Y\cap U_{2}}(U_{2},\alpha_{|U_{2}})$ and
$E=E'\oplus E'[1]$, so that $[E]=0$ in $K_{0,Y\cap
U_{2}}(U_{2},\alpha_{|U_{2}})$. Let $G=E_{|U_{12}}$. By the affine
part of the theorem (which we have already proven), we know that
there is $F\in Perf_{Y\cap U_{1}}(U_{1},\alpha_{|U_{1}})$ such
that $F_{|U_{12}}\simeq G$. Consider $\delta:Rj_{1*}F\oplus
Rj_{2*}E\longrightarrow Rj_{12*}G$ the adjunction morphism, and
let $D$ be the cocone of $\delta$, so that we have the
distinguished triangle
$$D\longrightarrow Rj_{1*}F\oplus Rj_{2*}E\stackrel{\delta}\longrightarrow
Rj_{12*}G\longrightarrow D[1].$$Applying to it the exact functor
$j_{1}^{*}$ we find that $j_{1}^{*}D\simeq F$, while applying
$j_{2}^{*}$ we get $j_{2}^{*}D\simeq E$. In this way we have shown
that $D\in Perf_{Y}(X,\alpha)$. By the hypothesis on $C$, we have:
$$0=Hom(D,C)=Hom(D,Rj_{2*}j_{2}^{*}C)=Hom(j_{2}^{*}D,j_{2}^{*}C)=Hom(E,j_{2}^{*}C)$$so
that $Hom(E',j_{2}^{*}C)=0$ for every $E'\in Perf_{Y\cap
U_{2}}(U_{2},\alpha_{|U_{2}})$. But $U_{2}$ verifies the theorem
by induction, so this implies $j_{2}^{*}C=0$, that is $C=0$.
\endproof

Using this theorem we can show the following:

\begin{cor}
\label{cor:class} Let $X$ be a smooth variety over a field $k$,
$\alpha\in BrX$, $Z$ a closed subset of $X$, $U=X\setminus Z$. Let
$\mathcal{F}\in Coh(U,\alpha_{|U})$. Then, there is
$\mathcal{E}\in Coh(X,\alpha)$ such that
$\mathcal{E}_{|U}\simeq\mathcal{F}$. In particular, $(X,\alpha)$
verifies the restriction condition.
\end{cor}

\proof Take $C=\mathcal{F}\oplus\mathcal{F}[1]$, so that $[C]=0$
in $K_{0}(U,\alpha_{|U})$. Since $X$ is smooth, $C\in
Perf(U,\alpha_{|U})$. By Theorem \ref{thm:obstruction}, we have
that there is $C'\in Perf(X,\alpha)=D^{b}(X,\alpha)$ such that
$C'_{|U}$ is quasi-isomorphic to $C$. Now, let
$\mathcal{E}=\mathcal{H}^{0}(C')\in Coh(X,\alpha)$. Since the
restriction is an exact functor, we have
$$\mathcal{E}_{|U}=\mathcal{H}^{0}(C'_{|U})\simeq
\mathcal{H}^{0}(C)=\mathcal{F}.$$\endproof

\section[Saturatedness of $D^{b}(X,\alpha)$]{Saturatedness of $D^{b}(X,\alpha)$}

Using the results of the previous section, following \cite{BdB} we
can even show that if $X$ is a smooth and proper variety over a
field $k$ and $\alpha\in BrX$, then $D^{b}(X,\alpha)$ is
saturated.

In \cite{BdB} is given the following definition

\begin{defn}
Let $\mathcal{T}$ a triangulated category such that for every
$A,B$ in $\mathcal{T}$ we have
$\sum_{i\in\mathbf{Z}}dim\,Hom(A,B[i])<\infty$. $\mathcal{T}$ is
called \textit{saturated} if every contravariant cohomological
functor of finite type $H:\mathcal{T}\longrightarrow Vect(k)$ is
representable.
\end{defn}

We shall need the following definitions and results. Here we use
the same notations as in section 3.1.

\begin{defn}
We say that a family of objects $\mathcal{E}\subseteq\mathcal{T}$
\textit{generates} (resp. \textit{strongly generates})
$\mathcal{T}$ if we have
$\langle\mathcal{E}\rangle_{\infty}=\mathcal{T}$ (resp. if there
is an integer $n$ such that
$\langle\mathcal{E}\rangle_{n}=\mathcal{T}$). Moreover, we say
that $\mathcal{T}$ is \textit{finitely generated} (resp.
\textit{finitely strongly generated}) if we can find a generating
family given by just one object, that will be called a
\textit{generator} (resp. \textit{strong generator}).
\end{defn}

\begin{prop}
\label{prop:karoubi} Let $\mathcal{T}$ be a triangulated category
which admits infinite direct sums, and let $\mathcal{T}^{c}$ the
full triangulated subcategory of compact objects. Then
$\mathcal{T}^{c}$ is Karoubian.
\end{prop}

\proof This is shown in \cite{BN}.\endproof

In \cite{BdB} it is show the following theorem:

\begin{thm}
\label{thm:saturated} Let $\mathcal{T}$ be a triangulated category
such that for every $A,B\in\mathcal{T}$ we have
$\sum_{i\in\mathbf{Z}}dim\,Hom(A,B[i])<\infty$. If $\mathcal{T}$
is Karoubian and is strongly finitely generated, then
$\mathcal{T}$ is saturated.
\end{thm}

\proof See proof of Theorem 1.3 in \cite{BdB}.\endproof

We want to use this criterion to show that if $X$ is a smooth
proper variety over a field $k$ and $\alpha\in BrX$, then
$D^{b}(X,\alpha)$ is saturated.

The fact that $X$ is proper implies that for any $A,B\in
D^{b}(X,\alpha)$ we have
$\sum_{i\in\mathbf{Z}}dim\,Hom(A,B[i])<\infty$. Moreover, the fact
that $X$ is smooth and that $\alpha\in BrX$ implies, by Theorem
\ref{thm:perf}, that $$D^{b}(X,\alpha)=Perf(X,\alpha)=
D(QCoh(X,\alpha))^{c},$$so that Proposition \ref{prop:karoubi}
tells us that $D^{b}(X,\alpha)$ is Karoubian. It remains to show
that we can find a strong generator. So, we start by proving that
we are able to find a generator for $Perf(X,\alpha)$.

\begin{prop}
\label{prop:generator} Let $X$ be a variety, $\alpha\in BrX$. Then
$D(QCoh(X,\alpha))$ is generated by a perfect complex. In
particular, $Perf(X,\alpha)$ is finitely generated.
\end{prop}

\proof By Lemma \ref{lem:neeman}, it suffices to show that there
is a perfect complex whose orthogonal is zero. The proof goes by
induction on the minimal number of affine open subschemes that
cover $X$. We start by $X=Spec\,A$. This was done in Theorem
\ref{thm:perf}: there, we showed that a generator is a locally
free $\alpha-$sheaf of finite rank (thought as a complex
concentrated in degree 0), that we will denote, as usual,
$\mathcal{E}$.

Now, let's suppose $X=U_{1}\cup U_{2}$, where $U_{1}=Spec\,A$ is
an open affine subscheme of $X$, and $U_{2}$ is an open subscheme
of $X$ which verifies the proposition. We will use the same
notations as in the proof of Theorem \ref{thm:obstruction}. In
particular, we are able to find a perfect generator $\mathcal{F}$
of $D(QCoh(U_{2},\alpha_{|U_{2}}))$. Now let
$\mathcal{F}'=\mathcal{F}\oplus\mathcal{F}[1]$, so that, from
Theorem \ref{thm:obstruction}, we know that there is a perfect
complex $P\in Perf(X,\alpha)$ such that
$j_{2}^{*}P\simeq\mathcal{F}'$.

Moreover, let $Y=X\setminus U_{2}=U_{1}\setminus U_{12}$, which is
a closed subscheme of $Spec\,A$, so that it will be given by
$f_{1},...,f_{r}\in A$. In the proof of Theorem
\ref{thm:obstruction} we showed that the complex $Q:=G_{r}^{\vee}$
associated to $Y$ is a perfect generator of
$D_{Y}(QCoh(U_{1},\alpha_{|U_{1}}))$.

We want to show is that $C=P\oplus Rj_{1*}Q$ is a perfect
generator of $D(QCoh(X,\alpha))$. By Theorem \ref{thm:perf} and
Lemma \ref{lem:neeman}, we know that it suffices to show that
$C^{\bot}=0$.

Since $Supp\,Q\subseteq Y$, we have that $j_{1}^{*}Rj_{1*}Q=Q$ and
$j_{2}^{*}Rj_{1*}Q=0$. In this way we see that $Rj_{1*}Q$ is in
$Perf(X,\alpha)$, so that $C\in Perf(X,\alpha)$. Moreover, using
Lemma \ref{lem:mayer} we can show that for every $D\in
D(QCoh(X,\alpha))$ we have
$$RHom(Rj_{1*}Q,D)=RHom(Q,j_{1}^{*}D).$$Now, suppose
$D\in C^{\bot}$. This gives $RHom(Rj_{1*}Q,D)=0$, so that
$j_{1}^{*}D\in Q^{\bot}$. As we saw in the proof of Theorem
\ref{thm:obstruction}, this implies that we get a canonical
isomorphism $D\simeq Rj_{2*}j_{2}^{*}D$.

But we have also that $RHom(P,D)=0$. This means
$$0=RHom(P,Rj_{2*}j_{2}^{*}D)=RHom(j_{2}^{*}P,j_{2}^{*}D)$$and
since $j_{2}^{*}P\simeq\mathcal{F}'$, this implies that
$j_{2}^{*}D$ is orthogonal to $\mathcal{F}$, which is a generator
of $D(QCoh(U_{2},\alpha_{|U_{2}}))$. This tells us that
$j_{2}^{*}D$ must be 0, and so $D=0$.\endproof

Now that we have shown that we can find a perfect generator for
the derived category $D(QCoh(X,\alpha))$, we can show the
following:

\begin{prop}
\label{prop:prodgen} Let $X,Y$ be two varieties, $\alpha\in BrX$,
$\beta\in BrY$, $\mathcal{F}$ a perfect generator of
$D(QCoh(X,\alpha))$, $\mathcal{G}$ a perfect generator of
$D(QCoh(Y,\beta))$. Then $\mathcal{F}\boxtimes\mathcal{G}$ is a
perfect generator of $D(QCoh(X\times Y,p^{*}\alpha\cdot
q^{*}\beta))$, where $p,q$ are the projections of $X\times Y$ on
$X$ and $Y$ respectively.
\end{prop}

\proof The fact that $\mathcal{F}\boxtimes\mathcal{G}$ is perfect
is clear. We have to show that if
$D\in\mathcal{F}\boxtimes\mathcal{G}^{\bot}=0$, then $D=0$. We
have, for every $i,j\in\mathbb{Z}$
$$0=Hom(p^{*}\mathcal{F}\otimes q^{*}\mathcal{G},D[i+j])=
Hom(p^{*}\mathcal{F},R\mathcal{H}om(q^{*}\mathcal{G},D[i])[j])=$$
$$=Hom(\mathcal{F},Rp_{*}R\mathcal{H}om(q^{*}\mathcal{G},D[i])[j]).$$Since
$\mathcal{F}$ is a generator for $D(QCoh(X,\alpha))$, we get
$Rp_{*}R\mathcal{H}om(q^{*}\mathcal{G},D[i])=0$ for every
$i\in\mathbb{Z}$. Now, take $U$ an open affine subscheme of $X$
and $\mathcal{E}$ a locally free $\alpha-$sheaf of finite rank in
$D(QCoh(X,\alpha))$. We have
$$Hom(\mathcal{E},Rp_{*}R\mathcal{H}om(q^{*}\mathcal{G},D[i]))=0$$so
that
$$0=Hom(\mathcal{E}_{|U},(Rp_{*}R\mathcal{H}om(q^{*}\mathcal{G},D[i]))_{|U}).$$Now,
let's denote $p_{U}$ and $q_{U}$ the projection of $U\times Y$ to
$U$ and $Y$. We get
$$0=Hom(p_{U}^{*}\mathcal{E},R\mathcal{H}om((q^{*}\mathcal{G})_{|U\times
Y},D[i]_{|U\times Y}))=$$
$$=Hom(p_{U}^{*}\mathcal{E}\otimes q_{U}^{*}\mathcal{G},D[i]_{|U\times
Y})=Hom(\mathcal{G},Rq_{U*}R\mathcal{H}om(p_{U}^{*}\mathcal{E},D[i]_{|U\times
Y})).$$Since $\mathcal{G}$ is a generator of $D(QCoh(Y,\beta))$
this implies
$$0=Rq_{U*}R\mathcal{H}om(p_{U}^{*}\mathcal{E},D[i]_{|U\times
Y})=Rq_{U*}R\mathcal{H}om(\mathcal{O}_{U\times
Y},p_{U}^{*}\mathcal{E}^{\vee}\otimes D[i]_{|U\times Y})$$the last
one being a sheaf on $Y$. Now, take $V$ an open affine subscheme
of $Y$, so that
$$Hom(\mathcal{O}_{V},Rq_{U*}R\mathcal{H}om(\mathcal{O}_{U\times
Y},(p^{*}\mathcal{E}^{\vee}\otimes D[i])_{|U\times Y}))=0$$that is
$$R\mathcal{H}om(\mathcal{O}_{X\times Y},p^{*}\mathcal{E}\otimes
D[i])_{|U\times V}=0$$for every $i\in\mathbb{Z}$, $U,V$ open
affine subschemes of $X$ and $Y$ respectively. We have found that
$R\mathcal{H}om(\mathcal{O}_{X\times
Y},p^{*}\mathcal{E}^{\vee}\otimes D[i])=0$, that implies
$p^{*}\mathcal{E}^{\vee}\otimes D=0$, and so $D=0$.\endproof

Now we can use Propositions \ref{prop:generator} and
\ref{prop:prodgen} to show the following:

\begin{prop}
\label{prop:strong} Let $X$ be a smooth variety over a field $k$,
$\alpha\in BrX$. Then $D^{b}(X,\alpha)$ has a strong generator.
\end{prop}

\proof Let $\gamma=p^{*}\alpha\cdot q^{*}\alpha^{-1}$. Since $X$
is smooth, we know that $X\times X$ is smooth and that the
structure sheaf $\mathcal{O}_{\Delta}$ of the diagonal is perfect.
We have that $\mathcal{O}_{\Delta}$ has structure of
$\gamma-$sheaf: if $\delta$ is the closed immersion of $\Delta$ in
$X\times X$, it's easy to show that $\delta^{*}\gamma=1$.

Let $\mathcal{F}$ be a perfect generator of $D(QCoh(X,\alpha))$,
and $\mathcal{G}$ a perfect generator of $D(QCoh(X,\alpha^{-1}))$
(we know that there are such generators from Proposition
\ref{prop:generator}). From Proposition \ref{prop:prodgen} we know
that $\mathcal{F}\boxtimes\mathcal{G}$ is a perfect generator of
$D(QCoh(X\times X,\gamma))$, so that there is $n\in\mathbb{N}$
such that $\mathcal{O}_{\Delta}\in
\langle\mathcal{F}\boxtimes\mathcal{G}\rangle_{n}$.

Now, we know that the Fourier-Mukai transform
$$\Phi_{\mathcal{O}_{\Delta}}:D^{b}(X,\alpha)\longrightarrow
D^{b}(X,\alpha)$$is the identity, so that for every $D\in
D(QCoh(X,\alpha))$, we have
$$D=Rp_{*}(q^{*}D\otimes\mathcal{O}_{\Delta})\in
\langle
Rp_{*}(q^{*}D\otimes(\mathcal{F}\boxtimes\mathcal{G}))\rangle_{n}.$$Since
$Rp_{*}(q^{*}D\otimes(\mathcal{F}\boxtimes\mathcal{G}))=\mathcal{F}\otimes
Rp_{*}q^{*}(D\otimes\mathcal{G})$, and $D\otimes\mathcal{G}$ is a
sheaf on $X$, we have that $$D\in \langle\mathcal{F}\otimes
R\Gamma(X,D\otimes\mathcal{G})\rangle_{n}=\overline{\langle\mathcal{F}\rangle_{n}}.$$Using
Lemma \ref{lem:neeman} we get
$D(QCoh(X,\alpha))^{c}=\langle\mathcal{F}\rangle_{n}$, and since
$X$ is smooth,
$D^{b}(X,\alpha)=\langle\mathcal{F}\rangle_{n}$.\endproof

\subsection*{Acknowledgment}
I mostly would like to thank Daniel Huybrechts and Raphael
Rouquier for encouragement and many helpful conversations during
the curse of this study, which is a continuation of my Mémoire de
DEA at the Université de Paris VII, during the year 2005. I would
also thank Paolo Stellari and Alberto Canonaco for having read
preliminary versions of this work, and having signaled me mistakes
and improvements.

\noindent Arvid Perego, Laboratoire de Mathématiques Jean Leray,
Université de Nantes, 2, rue de la Houssinière, BP 92208, F-44322
Nantes Cedex 03, France

\noindent \textit{E-mail address}:
arvid.perego@math.univ-nantes.fr

\end{document}